\documentclass[10pt,twocolumn,twoside]{IEEEtran}

\usepackage{amsmath}
\usepackage{amssymb}
\usepackage{amsthm}
\usepackage{algorithm}
\usepackage{algorithmic}
\usepackage{tikz}
\usetikzlibrary{calc,decorations.pathmorphing}
\usepackage{flushend}

\theoremstyle{plain}
\newtheorem{thm}{Theorem}

\theoremstyle{definition}
\newtheorem{rem}{Remark}

\newcommand{\bN}{{\mathbb{N}}}

\newcommand{\cE}{{\mathcal{E}}}
\newcommand{\cG}{{\mathcal{G}}}
\newcommand{\cV}{{\mathcal{V}}}

\newcommand{\Pre}{\mathrm{Pre}}
\newcommand{\Post}{\mathrm{Post}}

\title{Finite time distributed averaging over ring networks}

\author{Alessandro Falsone, Kostas Margellos, Simone Garatti, Maria Prandini%
	\thanks{Corresponding author: Alessandro Falsone.
	}%
	\thanks{Alessandro Falsone, Simone Garatti, and Maria Prandini are with Dipartimento di Elettronica, Informazione e Bioingegneria, Politecnico di Milano, Milano, Italy (phone: +39 02 2399 \{4028,3650,3441\} e-mail: \texttt{\{alessandro.falsone, simone.garatti, maria.prandini\}@polimi.it}).
	}%
	\thanks{Kostas Margellos is with Department of Engineering Science, University of Oxford, Oxford, United Kingdom (phone: +44 1865 283912; e-mail: \texttt{kostas.margellos@eng.ox.ac.uk).}
	}%
	\thanks{This work is partially supported by the European Commission under the project UnCoVerCPS with grant number 643921.
	}%
}%

\begin{document}

\maketitle
\thispagestyle{empty}
\pagestyle{empty}

\begin{abstract}
We consider a multi-agent system where each agent has its own estimate of a given quantity and the goal is to reach consensus on the average. To this purpose, we propose a distributed consensus algorithm that guarantees convergence to  the average in a finite number of iterations. The algorithm is tailored to ring networks with bidirectional pairwise communications. If the number of agents $m$ is even, say $m=2n$, then, the number of iterations needed is equal to $n$, which in this case is the diameter of the network, whereas the number of iterations grows to $3n$ if the number of agents is odd and equal to $m=2n+1$.
\end{abstract}


\section{Introduction} \label{sec:intro}
A typical problem encountered in a multi-agent system is that all agents are aiming at some common goal but they can communicate only with their neighbors to this purpose.
Achieving a common goal often translates into reaching an agreement on the value taken by some quantity (the agreement variable) via some consensus algorithm.  Basic consensus algorithms date back to \cite{tsitsiklis1984problems}. Only in the last two decades,  they have attracted the attention of both the computer science \cite{lynch1996distributed} and control engineering \cite{morse2003coordination} communities.

There are different forms of consensus. Here, we are concerned with distributed averaging where each agent has its own estimate of a certain quantity and the goal is to reach consensus on the average of the values stored by all agents. The two most common approaches in the literature to the distributed averaging  problem are linear iterations and gossip algorithms. The former is an iterative scheme where each agent updates its estimate of the average by taking a linear combination of its current estimate and those received by its neighboring agents, \cite{boyd2004fast}. Gossip algorithms are similar, but they require that only pairwise communication occurs, \cite{boyd2005gossip, morse2011deterministic}. In both approaches, the interactions between agents can be modeled as a time-varying weighted graph, with vertices and edges respectively representing the agents and the communication links, and weights assigned to edges being the coefficients of the linear combinations. Under mild assumptions on the coefficients and on the structure of the graph across iterations, both linear iterations schemes and gossip algorithms have been proven to asymptotically converge to the average, \cite{boyd2004fast}. Since reaching consensus asymptotically can be limiting in practice, there is a fair amount of literature on how to analyze and optimize the convergence rate of both algorithms, see \cite{boyd2004fast, olshevsky2009convergence} and references therein.

As pointed out in \cite{lynch1996distributed, boyd2004fast, morse2011deterministic}, distributed average computation can be actually solved in a finite number of iterations by making each agent keep collecting all the values received and passing them to its neighbors. After a number of iterations equal to the diameter of the graph (i.e., the maximum distance between any two vertices), every agent knows all the values and can compute their average. This algorithm is apparently very simple but it actually presents two main issues: i) the number of values that each agent needs to store (and thus the memory usage) grows linearly with the number of agents, and ii) an unnecessary amount of information is exchanged which might overload the communication channels. If the quantities to be averaged were vectors instead of scalars,  the previous two issues would become even more critical.
Much effort has then be devoted in the literature to design simple algorithms (like linear iterations) that are able to solve the distributed averaging problem in finite time, see \cite{sundaram2007finite, ko2009matrix, georgopoulos2011definitive, kibangou2011finite, hendrickx2014graph, morse2014finite, hendrickx2015finite, oliva2016distributed} just to name a few.
Most of these works take a global perspective, in that they assume to know the topology of the (undirected) graph, and then orchestrate the weights so as to achieve finite time convergence.
In \cite{sundaram2007finite} the topology is assumed to be time-invariant, and the weights are kept constant across iterations. The solution to the finite time averaging problem is then given in terms of the minimal polynomial of some matrix which gathers the weights and matches the graph topology. Convergence is  achieved in $D+1$ steps, where $D$ is the degree of the minimal polynomial, provided that each agent stores all the $D+1$ previous values, which, however, might be impractical for large networks. A further result in \cite{sundaram2007finite} is that the degree of the minimal polynomial can be computed in a distributed way in a number of iterations that is equal to the number of agents.
In \cite{kibangou2011finite} the network is also assumed to be time-invariant, but the weights (and thus the weight matrix) change across iterations. An analytic solution based on the joint diagonalization of the weight matrices is provided, and it is shown that convergence can be achieved in a number of steps equal to the number of distinct eigenvalues of the Laplacian matrix of the graph.

In \cite{oliva2016distributed} the authors consider a fixed network topology and develop an algorithm based on a two-stage max-consensus which achieves finite time convergence in $d(2m+1)$ iterations, where $d$ is the diameter of the graph and $m$ is the number of agents.
In \cite{ko2009matrix, morse2014finite, hendrickx2015finite}, the authors focus only on trees, i.e., graphs without loops, and show how to reach finite time convergence over these topologies. Note that, even though it is always possible to construct a spanning tree (i.e., a tree which reaches all vertices) from a connected graph by dropping some links, the diameter of the resulting tree might be greater than the diameter of the original graph, thus requiring more iterations to reach convergence.
In \cite{georgopoulos2011definitive} the finite time averaging problem is rephrased as a nonlinear optimization program that can be solved numerically, where the weights in the graph are optimization variables. Finite time average can be achieved in this case with a number of iterations that ranges between $d$ and $2d$, where $d$ is the diameter of the network. However, agents have to solve a nonlinear program in a distributed fashion, which might be difficult from a computational point of view.
In \cite{hendrickx2014graph} the authors prove that for specific type of graphs finite time averaging can be achieved in a number of steps equal to the diameter of the graph. Ring networks are included in this setting, but the communication graph is kept constant and gossip constraints are hence not considered.
Finally, in \cite{ko2009matrix,shi2012finite}, finite time averaging under the gossip constraint has been investigated. In particular, in \cite{shi2012finite} it is proven that if the weights are constant and equal to $1/2$, then finite time convergence on undirected graphs
can be achieved only if the number of vertices equals a power of $2$. Finite time averaging can be achieved on a directed graph with weights equal to $1/2$, but in such a case the number of iteration required for a network with $m=2^n+r$ agents, where $r\in\bN$, is $mn+2r$, which grows more than linearly with $m$.
In \cite{ko2009matrix}, it is shown that finite time averaging can be achieved under the gossip constraint by allowing for time-varying coefficients. An upper bound that scales quadratically with the number $m$ of agents is given for the number of iterations required to achieve finite time averaging.
% but their algorithms are based on the construction of a spanning tree.

In this work we are concerned with the design of a finite time averaging algorithm specifically tailored to a ring topology, subject to a gossip constraint. By allowing the weights of the communication graph to be time-varying, we are able to prove finite time convergence for a network with an arbitrary number of agents.
More specifically, for ring networks with an even number of agents, we propose an  algorithm that ensures finite time convergence and attains the lower bound on the number of iterations needed for a general communication graph with no gossip constraint, which is the diameter of the graph. As for ring networks with an odd number of agents, an algorithm based on the case of an even number of vertices is designed. Still convergence in finite time to the average is guaranteed in a number of iterations that scales linearly with the number of agents. In this case, however, the lower bound is not attained. As for memory requirements, each agent needs to store only its estimate of the mean when the total number of agents is even, whereas the memory requirement is doubled in the odd case. In both cases, the required memory is independent of the number of agents.

%The gossip constraint requires the network to be time-varying, and to the best of the authors' knowledge finite time gossiping on a time-varying ring topology has never been addressed before. The contributions of the paper are next stated. For rings with an even number of vertices our algorithm ensure finite time convergence and attain the lower bound on the number of iteration, which is the diameter of the ring. As for rings with an odd number of vertices we propose a finite time algorithm based on the even ring case, but this time the lower bound is not attained. \AF{I am not sure about this last statement, because gossiping adds constraints on the communication structure, and thus on the minimum number of iterations needed.}

The rest of the paper is organized as follows. In Section \ref{sec:setup}, we introduce the problem set-up with the ring topology and the gossip constraint. The proposed distributed averaging algorithms for the cases of even and odd number of agents are described in Section \ref{sec:solution}, where finite time convergence and related bounds on the number of iterations are shown. Finally, in Section \ref{sec:conclusions} some concluding remarks are drawn.

\section{Problem setup} \label{sec:setup}

We address the case of a multi-agent system characterized by a ring communication network, where each agent can communicate only with one of each neighbors at a time (gossip constraint) to take some joint decision. The agents are $m$ in total and each one is identified by an integer $i$ taking values in $\{1, 2,\dots,m\}$.  Each agent $i$ has its own estimate $x_i(0)$ of some quantity of interest. The goal is to devise a distributed algorithm through which each agent is able to compute the average
\begin{align}\label{eq:x_bar}
	\bar{x} = \frac{1}{m}\sum_{i=1}^{m} x_i(0),
\end{align}
in a finite number of iterations only by exchanging information with its neighbors under the gossip constraint. Evidently, only the case with $m>2$ is of interest.

The network communication structure can be represented as an undirected graph $\cG = (\cV,\cE)$, where $\cV = \{1,\dots,m\}$ is the set of vertices representing the agents, and $\cE$ is the set of edges, defined as the following collection of ordered pairs of vertices:
\begin{equation*}
	\cE = \bigcup_{i=1}^{m-1} \{ (i,i+1), (i+1,i) \} \cup \{ (1,m),(m,1) \},
\end{equation*}
according to the ring topology. Since we are dealing with an undirected graph, with a slight abuse of notation we will use $(i,j)$ to denote both edges $(i,j)$ and $(j,i)$. According to this notation, for $m>2$, the number of edges is even if $m$ is even, odd otherwise.

We next define the communication protocol used by the agents to exchange information over the network.
We are here concerned with the design of a synchronous pairwise communication strategy. The agents communicate in rounds, and successive communication rounds are indexed by an integer $k\in\bN$, which is also referred to as time step or iteration. At any round $k$, each agent is allowed to communicate only with one of its neighbors (gossip constraint). Since communication channels are represented by the edges of $\cG$, the design of the communication protocol then reduces to specify a proper sequence of edges to activate, where edge $(i,j)$ is said to be \textit{active} at time step $k$ if agents $i$ and $j$ exchange information at time $k$. Of course, more edges can be activate at the same iteration, as long as they do not have any vertex in common, so as to comply with the pairwise communication constraint (see the multigossip framework in \cite{morse2014finite}).

\begin{figure}[t]
	\centering
	\begin{tikzpicture}
		\def\RAD{0.35}
		\def\DIST{2}
		\def\ANGLE{360/6}
		
		% Define colors
		\definecolor{blue}{rgb}{0.0000,0.4470,0.7410}
		\definecolor{red}{rgb}{0.9500,0.3250,0.0980}
		
		\coordinate (O) at (0,0);
		
		\coordinate (N1) at ($(O)+(90+0*\ANGLE:\DIST)$);
		\coordinate (N2) at ($(O)+(90+1*\ANGLE:\DIST)$);
		\coordinate (N3) at ($(O)+(90+2*\ANGLE:\DIST)$);
		\coordinate (N4) at ($(O)+(90+3*\ANGLE:\DIST)$);
		\coordinate (N5) at ($(O)+(90+4*\ANGLE:\DIST)$);
		\coordinate (N6) at ($(O)+(90+5*\ANGLE:\DIST)$);
		
		\draw[thick,blue] (N1) -- (N2);
		\draw[thick,red,style={decorate, decoration=snake}] (N2) -- (N3);
		\draw[thick,blue] (N3) -- (N4);
		\draw[thick,red,style={decorate, decoration=snake}] (N4) -- ($0.6*(N4)+0.4*(N5)$);
		\draw[thick,red,dotted] ($0.6*(N4)+0.4*(N5)$) -- ($0.4*(N4)+0.6*(N5)$);
		\draw[thick,blue] ($0.4*(N5)+0.6*(N6)$) -- (N6);
		\draw[thick,blue,dotted] ($0.4*(N5)+0.6*(N6)$) -- ($0.6*(N5)+0.4*(N6)$);
		\draw[thick,red,style={decorate, decoration=snake}] (N6) -- (N1);
		
		\filldraw[fill=white] (N1) circle (\RAD) node {$1$};
		\filldraw[fill=white] (N2) circle (\RAD) node {$2$};
		\filldraw[fill=white] (N3) circle (\RAD) node {$3$};
		\filldraw[fill=white] (N4) circle (\RAD) node {$4$};
		\node at (N5) {$\dots$};
		\filldraw[fill=white] (N6) circle (\RAD) node {$m$};
	\end{tikzpicture}
	\caption{Communication protocol under the gossip constraint in the case of an even number of agents: edges are grouped in two sets by a color coding. Only edges with the same color are active simultaneously. }
	\label{fig:communication_protocol_even}
\end{figure}
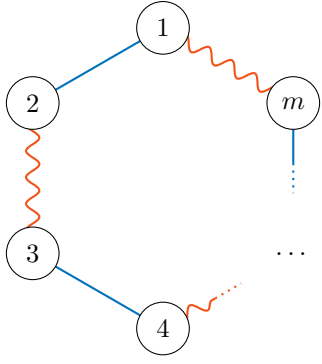

To reduce the number of communication rounds we need to parallelize as much as possible the number of simultaneous pairwise communications. This can be interpreted as an edge-coloring problem on $\cG$, \cite{morse2014finite}, where edges with the same color represent communication channels which can be active at the same time. For a ring communication network with an even number $m>2$ of vertices  we need just two colors to minimize the number of communication rounds needed for all edges to be activated while satisfying the gossip constraint. In the odd case we need at least three colors. Figure~\ref{fig:communication_protocol_even} shows  the coloring scheme for the even case in a pictorial form, with two groups of edges characterized by two different colors, blue and red, that can be activated alternatively, e.g., all the blue straight edges are activate at iterations $k=1,3,5, \dots$, whereas the red wavy ones at $k=2,4,6\dots$.

Finally, we need to specify the weights associated with the edges, which define the linear update of the agents local estimate of the mean based on the information received from its neighbors. To this end, assume that at time $k$, agents $i$ and $j$ communicate. Agent $i$ performs then the following update for its local estimate of the mean:
\begin{align}
		&x_i(k) = (1-\alpha_k) x_i(k-1) + \alpha_k x_j(k-1), \label{eq:update_step}	
\end{align}
where $\alpha_k$ is the (time-varying) weight associated with edge $(i,j)$ at time $k$ entering the convex update rule. Agent $j$ performs an analogous update step.

The problem to be addressed is how to appropriately select $\alpha_k$, $k=1,2, \dots$, so as to guarantee that there exists a finite time $T$ such that $x_i(T) = \bar{x}$ for all $i=1,\dots,m$, where $\bar{x}$ is given by \eqref{eq:x_bar}.

\section{Proposed solution} \label{sec:solution}

We next provide a solution to the finite time distributed averaging problem over a ring network topology with pairwise communication constraints as described in the previous section.
In particular, we shall define which edges are activated at each step $k$ and appropriately set the values for the $\alpha_k$ weights in the update rule \eqref{eq:update_step} to achieve finite time convergence of all agents' estimates to the average \eqref{eq:x_bar}.

We will first address the case when the number of agents in the ring is even, and then we will extend the result to rings with an odd number of agents.

\subsection{Ring with an even number of agents} \label{sec:solution_even}
For a ring with an even number of agents ($m = 2n$), the coefficients of the linear combination can be chosen according to the following rule:
\begin{equation}
	\alpha_k =
	\begin{cases}
		\frac{k}{k+1}, & 1 \le k < n \\
		\frac{1}{2}, & k = n \\
        0, & k > n
	\end{cases}
	\label{eq:weighting_coeff}
\end{equation}
which due to \eqref{eq:update_step}	entails that the estimate is kept constant after $k=n$ steps and hence the iterative process can be terminated.

The distributed overall procedure is obtained by making each agent run synchronously with the others Algorithm~\ref{algo:agent_algorithm}, which involves alternatively communicating with the predecessor and successor neighbor in the ring topology (according to the agent indexing shown in Figure \ref{fig:communication_protocol_even}), while updating its estimate  according to \eqref{eq:update_step} with the coefficients of the convex combination defined as in \eqref{eq:weighting_coeff}. Agents should agree on the neighbor to communicate with at the first iteration, i.e., which of the two edge groups of Figure~\ref{fig:communication_protocol_even} should be simultaneously activated.
This choice is actually embedded in Algorithm~\ref{algo:agent_algorithm} (see line 4) and it is simply based on the agent identification number $i$.
Functions $\Post, \Pre$: $\{1,2, \dots, 2n\} \to \{1,2, \dots, 2n\}$ in Algorithm~\ref{algo:agent_algorithm} serve the purpose of specifying the successor and predecessor neighbor, respectively, and are given by:
\begin{equation}
	\Post(i) =
	\begin{cases}
		i+1,	& i = 1,\dots,2n-1 \\
		1,		& i = 2n
	\end{cases}
	\label{eq:right_function}
\end{equation}
and
\begin{equation}
	\Pre(i) =
	\begin{cases}
		2n,		& i = 1\\
        i-1,		& i = 2,\dots,2n.
	\end{cases}
	\label{eq:left_function}
\end{equation}
According to Algorithm~\ref{algo:agent_algorithm}, at iteration $k=1$ agent $i=1$ communicates with its successor agent $2$, and, hence, the blue colored edges in Figure \ref{fig:communication_protocol_even} are activated.\\
The following theorem holds for the proposed distributed scheme.

\begin{algorithm}[t]
	\caption{Algorithm for agent $i$ -- $m$ even}
	\begin{algorithmic}[1]
		\STATE $x_i(0) \gets$ initial value for agent $i$
		\FOR{$k=1$ \TO $n$}
			\STATE \textit{\% Select a neighbor to communicate with}
			\IF{$i+k$ is even}
				\STATE $j \gets \Post(i)$
			\ELSE
				\STATE $j \gets \Pre(i)$
			\ENDIF
			\STATE \textit{\% Update the estimate using \eqref{eq:update_step} and \eqref{eq:weighting_coeff}}
			\STATE $x_i(k) \gets (1-\alpha_k) x_i(k-1) + \alpha_k x_j(k-1)$
		\ENDFOR
		\RETURN $x_i(n)$
	\end{algorithmic}
	\label{algo:agent_algorithm}
\end{algorithm}

\begin{thm}[Finite time consensus]\label{thm:correctness_algorithm}
Given a ring network with $m=2n$ agents, if all of them apply Algorithm~\ref{algo:agent_algorithm} synchronously, then, after $n$ iterations, the estimate $x_i$ computed by each agent $i$, $i=1,2,\dots,2n$, equals the average of their initial values \eqref{eq:x_bar}.

\begin{proof}
Our aim is to show that by applying Algorithm~\ref{algo:agent_algorithm} to a ring with $2n$ agents we get
\begin{equation}
	x_i(n) = \frac{1}{2n}\sum_{j=1}^{2n} x_j(0),
	\label{eq:average}
\end{equation}
for every  agent $i$, $i=1,2,\dots,2n$.

%To this purpose we need to introduce some notation. Let us consider the agent that we shall reach by moving $h$-th positions to the right of agent $i$. Formally, this is agent  $\RIGHT^h(i)$ whose index is obtained by applying the $\RIGHT(\cdot)$ function in \eqref{eq:right_function} $h$ times. Correspondingly, $x_{\RIGHT^h(i)}$ is its estimate of the average.  To ease the notation, we will replace $\RIGHT^h(i)$ by $i+h$ and  $x_{\RIGHT^h(i)}$ by $x_{i+h}$. Similarly, we shall refer to the $h$-th agent on the left of agent $i$ and  its estimate as $i-h$ and $x_{i-h}$, in place of $\LEFT^h(i)$ and $x_{\LEFT^h(i)}$.

It suffices to show that \eqref{eq:average} holds for $i$ such that $i+n$ is even (i.e., $i$ even if $n$ is even, odd if $n$ is odd). This is due to the fact that Algorithm~\ref{algo:agent_algorithm} agent $i+1$ (for which $i+1+n$ is odd) communicates with agent $i$ at step $k=n$ and its estimate satisfies
\begin{align*}
	x_{i+1}(n)  &= (1-\alpha_n) x_{i+1}(n-1) + \alpha_n x_i(n-1),
\end{align*}
which is equal to $x_i(n) =(1-\alpha_n) x_{i}(n-1) + \alpha_n x_{i+1}(n-1)$ given that $\alpha_n= \frac{1}{2}$. Then, all pairs of agents $(i,i+1)$ will converge to the average in $n$ steps.

We shall then focus next on agent $i$ such that $i+n$ is even, and prove that the following equation hold
\begin{align}
		x_i(n) &= \frac{1}{2n} \sum_{j = i-s+1}^{i+s} x_j(n-s-1) \nonumber \\
			   &+ \frac{n-s}{2n} \Big[ x_{i-s}(n-s-1) + x_{i+s+1}(n-s-1) \Big], \label{eq:generic_step}
\end{align}
for all $s = 1,\dots,n-1$. Note that, by substituting $s = n-1$ in \eqref{eq:generic_step}, we get \eqref{eq:average}, which concludes the proof of the theorem.

The proof of equation \eqref{eq:generic_step} is by induction, i.e., we show that it is true for $s=1$ (step 1), assume that it holds for some $s$ (induction hypothesis), and then show that this is also the case for $s+1$ (step 2).

\vspace*{1mm}

{\bf Step 1:}
By Algorithm~\ref{algo:agent_algorithm},
\begin{equation}
	x_i(n) = \frac{1}{2} x_i(n-1) + \frac{1}{2} x_{i+1}(n-1).
	\label{eq:last_step}
\end{equation}
At step $k=n-1$, we have that $i+k = i+n-1$ is odd and $i+1+k = i+n$ is even, hence
\begin{align}
&	x_i(n-1) = \frac{1}{n} x_i(n-2) + \frac{n-1}{n} x_{i-1}(n-2),
	\label{eq:pre_last_step_i}\\
&	x_{i+1}(n-1) = \frac{1}{n} x_{i+1}(n-2) + \frac{n-1}{n} x_{i+2}(n-2),
	\label{eq:pre_last_step_ip1}
\end{align}
and, by substituting \eqref{eq:pre_last_step_i} and \eqref{eq:pre_last_step_ip1} into \eqref{eq:last_step}, we obtain
	\begin{align*}
		x_i(n) = &\frac{1}{2n} \Big[ x_i(n-2) + x_{i+1}(n-2) \Big] \nonumber \\
			   &+ \frac{n-1}{2n} \Big[ x_{i-1}(n-2) + x_{i+2}(n-2) \Big].
	\end{align*}
which is \eqref{eq:generic_step} with $s = 1$.

\smallskip

{\bf Step 2:} For the sake of clarity we will treat the two terms in \eqref{eq:generic_step} separately: the summation first and then the other term.

\smallskip

\emph{First term:}
Let us start by considering the first two contributions in the summation, that is:
\begin{equation}
	x_{i-s+1}(n-s-1) + x_{i-s+2}(n-s-1).
	\label{eq:generic_step_summatio_two_terms}
\end{equation}
By Algorithm~\ref{algo:agent_algorithm} when $k = n-s-1$ (and, hence, $i-s+1+k = i+n-2s$ is even and $i-s+2+k = i+n-2s+1$ is odd), we get
	\begin{align}
		x_{i-s+1}(n-s-1) &= \frac{1}{n-s} x_{i-s+1}(n-s-2) \nonumber \\
						 &+ \frac{n-s-1}{n-s} x_{i-s+2}(n-s-2),\label{eq:generic_step_summation_1st_term}
	\end{align}

	\begin{align}
		x_{i-s+2}(n-s-1) &= \frac{1}{n-s} x_{i-s+2}(n-s-2) \nonumber \\
						 &+ \frac{n-s-1}{n-s} x_{i-s+1}(n-s-2).\label{eq:generic_step_summation_2nd_term}
	\end{align}

Substituting \eqref{eq:generic_step_summation_1st_term} and \eqref{eq:generic_step_summation_2nd_term} into \eqref{eq:generic_step_summatio_two_terms} we get that
	\begin{align}
		x_{i-s+1}&(n-s-1) + x_{i-s+2}(n-s-1) \nonumber \\
						 &= x_{i-s+1}(n-s-2) + x_{i-s+2}(n-s-2).\label{eq:generic_step_two_term_equivalence}
	\end{align}

The time-invariance property for the pairwise summation in \eqref{eq:generic_step_two_term_equivalence} holds true for all other pairs in the summation of the first term in \eqref{eq:generic_step}, which always contains an even number of terms, i.e., $2s$. We can thus conclude that
\begin{equation}
	\sum_{j = i-s+1}^{i+s} x_j(n-s-1) = \sum_{j = i-s+1}^{i+s} x_j(n-s-2).
	\label{eq:generic_step_summation_equivalence}
\end{equation}

\emph{Second term:}
Now consider the second term in \eqref{eq:generic_step}, recalled here for ease of reference:
\begin{equation}
	\frac{n-s}{2n} \Big[ x_{i-s}(n-s-1) + x_{i+s+1}(n-s-1) \Big].
	\label{eq:generic_step_2nd_term}
\end{equation}
 By Algorithm~\ref{algo:agent_algorithm} when $k = n-s-1$ (and, hence, $i-s+k = i+n-2s-1$ is odd, while $i-s+1+k = i+n-2s$ is even), we get
	\begin{align}
		x_{i-s}(n-s-1) =& \frac{1}{n-s} x_{i-s}(n-s-2) \nonumber \\
					   &+ \frac{n-s-1}{n-s} x_{i-s-1}(n-s-2),\label{eq:generic_step_2nd_term_1st_part}
	\end{align}

	\begin{align}
		x_{i+s+1}(n-s-1) =& \frac{1}{n-s} x_{i+s+1}(n-s-2) \nonumber \\
						 &+ \frac{n-s-1}{n-s} x_{i+s+2}(n-s-2).\label{eq:generic_step_2nd_term_2nd_part}
	\end{align}
 Substituting \eqref{eq:generic_step_2nd_term_1st_part} and \eqref{eq:generic_step_2nd_term_2nd_part} into \eqref{eq:generic_step_2nd_term}, \eqref{eq:generic_step_2nd_term} reduces to
	\begin{align}
%		&\frac{n-s}{2n} \Big[ x_{i-s}(n-s-1) + x_{i+s+1}(n-s-1) \Big] = \\
					&\frac{1}{2n} \Big[ x_{i-s}(n-s-2) + x_{i+s+1}(n-s-2) \Big] \nonumber \\
					&+ \frac{n-s-1}{2n} \Big[ x_{i-s-1}(n-s-2) + x_{i+s+2}(n-s-2) \Big].\label{eq:generic_step_2nd_term_equivalence}
	\end{align}

Finally, substituting \eqref{eq:generic_step_summation_equivalence} and \eqref{eq:generic_step_2nd_term_equivalence} into \eqref{eq:generic_step}, , which we assume that holds due to our induction hypothesis, we have that
	\begin{align*}
		&x_i(n) = \frac{1}{2n} \sum_{j = i-s}^{i+s+1} x_j(n-s-2) \nonumber \\
			   &+ \frac{n-s-1}{2n} \Big[ x_{i-s-1}(n-s-2) + x_{i+s+2}(n-s-2) \Big],
	\end{align*}
which is the same expression as \eqref{eq:generic_step} with $s+1$ in place of $s$, thus concluding the proof by induction of \eqref{eq:generic_step}.

\end{proof}
\end{thm}

\begin{rem}[Speed of convergence]
	Since the diameter  of a ring network with $2n$ agents is $n$, we need at least $n$ iterations to make sure that each agent receives information from all the other agents.
    This implies that the result of Theorem \ref{thm:correctness_algorithm}, which shows that we need exactly n iterations to converge to the average in ring networks with an even number of agents, is tight. To the best of our knowledge, this outperforms existing results in the literature.
%Thus, we can conclude that there does not exist any algorithm which is faster than Algorithm~\ref{algo:agent_algorithm} in solving the distributed averaging problem over ring networks with an even number of agents.
\end{rem}

\begin{rem}[Efficiency]
Since at each iteration every agent communicates with a single agent only, the information transmitted in the proposed distributed algorithm is limited compared to alternative solutions in the literature with either fixed or time-varying topology without gossip constraints. %with fixed topology. E QUELLE CON TOPOLOGIA VARIABILE?
\end{rem}

\begin{rem}[Interpretation]
At time $k$, the update of agent $i$, $i=1,\ldots,m$, in \eqref{eq:update_step} can be alternatively written as
\begin{align}
x_i(k) = x_i(k-1) -\alpha_k (x_i(k-1) - x_j(k-1)). \label{eq:update_step_alter}	
\end{align}
By inspection of \eqref{eq:update_step_alter}, it can be observed that the estimate of the average that each agent maintains evolves as a discrete time integrator, where the quantity that gets integrated is the mismatch between the communicating agents' estimates, weighted by $\alpha_k$. This can be thought of as the evolution of a closed-loop dynamical system, with the feedback gain fixed according to \eqref{eq:weighting_coeff}. The interpretation of Theorem \ref{thm:correctness_algorithm}, is that the dynamical system \eqref{eq:update_step_alter}, reaches $\bar{x}$ at the $n$-th step, for all $i=1,\ldots,m$.
\end{rem}

\subsection{Ring with an odd number of agents} \label{sec:solution_odd}

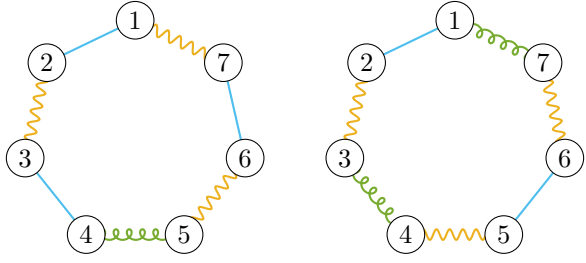
\begin{figure}[t]
	\centering
	\begin{tikzpicture}
		\def\RAD{0.25}
		\def\DIST{1.5}
		\def\ANGLE{360/7}
		
		% Define colors
		\definecolor{blue}{rgb}{0.3010,0.7450,0.9330}
		\definecolor{red}{rgb}{0.9290,0.6940,0.1250}
		\definecolor{green}{rgb}{0.4660,0.6740,0.1880}

		% First configuration
		\coordinate (O) at (0,0);
		
		\coordinate (N1) at ($(O)+(90+0*\ANGLE:\DIST)$);
		\coordinate (N2) at ($(O)+(90+1*\ANGLE:\DIST)$);
		\coordinate (N3) at ($(O)+(90+2*\ANGLE:\DIST)$);
		\coordinate (N4) at ($(O)+(90+3*\ANGLE:\DIST)$);
		\coordinate (N5) at ($(O)+(90+4*\ANGLE:\DIST)$);
		\coordinate (N6) at ($(O)+(90+5*\ANGLE:\DIST)$);
		\coordinate (N7) at ($(O)+(90+6*\ANGLE:\DIST)$);
		
		\draw[thick,blue] (N1) -- (N2);
		\draw[thick,red,style={decorate, decoration={snake, segment length=5, amplitude=2}}] (N2) -- (N3);
		\draw[thick,blue] (N3) -- (N4);
		\draw[thick,green,style={decorate, decoration={coil, aspect=1, segment length=5, amplitude=2}}] (N5) -- (N4);
%		\draw[thick,green!75!black,dashed] (N4) -- (N5);
		\draw[thick,red,style={decorate, decoration={snake, segment length=5, amplitude=2}}] (N5) -- (N6);
		\draw[thick,blue] (N6) -- (N7);
		\draw[thick,red,style={decorate, decoration={snake, segment length=5, amplitude=2}}] (N7) -- (N1);
		
		\filldraw[fill=white] (N1) circle (\RAD) node {$1$};
		\filldraw[fill=white] (N2) circle (\RAD) node {$2$};
		\filldraw[fill=white] (N3) circle (\RAD) node {$3$};
		\filldraw[fill=white] (N4) circle (\RAD) node {$4$};
		\filldraw[fill=white] (N5) circle (\RAD) node {$5$};
		\filldraw[fill=white] (N6) circle (\RAD) node {$6$};
		\filldraw[fill=white] (N7) circle (\RAD) node {$7$};
		
		% Second configuration
		\coordinate (O_2) at (4.25,0);
		
		\coordinate (N1_2) at ($(O_2)+(90+0*\ANGLE:\DIST)$);
		\coordinate (N2_2) at ($(O_2)+(90+1*\ANGLE:\DIST)$);
		\coordinate (N3_2) at ($(O_2)+(90+2*\ANGLE:\DIST)$);
		\coordinate (N4_2) at ($(O_2)+(90+3*\ANGLE:\DIST)$);
		\coordinate (N5_2) at ($(O_2)+(90+4*\ANGLE:\DIST)$);
		\coordinate (N6_2) at ($(O_2)+(90+5*\ANGLE:\DIST)$);
		\coordinate (N7_2) at ($(O_2)+(90+6*\ANGLE:\DIST)$);
		
		\draw[thick,blue] (N1_2) -- (N2_2);
		\draw[thick,red,style={decorate, decoration={snake, segment length=5, amplitude=2}}] (N2_2) -- (N3_2);
		\draw[thick,green,style={decorate, decoration={coil, aspect=1, segment length=5, amplitude=2}}] (N4_2) -- (N3_2);
%		\draw[thick,green!75!black,dashed] (N3_2) -- (N4_2);
		\draw[thick,red,style={decorate, decoration={snake, segment length=5, amplitude=2}}] (N4_2) -- (N5_2);
		\draw[thick,blue] (N5_2) -- (N6_2);
		\draw[thick,red,style={decorate, decoration={snake, segment length=5, amplitude=2}}] (N6_2) -- (N7_2);
		\draw[thick,green,style={decorate, decoration={coil, aspect=1, segment length=5, amplitude=2}}] (N1_2) -- (N7_2);
%		\draw[thick,green!75!black,dashed] (N7_2) -- (N1_2);
		
		\filldraw[fill=white] (N1_2) circle (\RAD) node {$1$};
		\filldraw[fill=white] (N2_2) circle (\RAD) node {$2$};
		\filldraw[fill=white] (N3_2) circle (\RAD) node {$3$};
		\filldraw[fill=white] (N4_2) circle (\RAD) node {$4$};
		\filldraw[fill=white] (N5_2) circle (\RAD) node {$5$};
		\filldraw[fill=white] (N6_2) circle (\RAD) node {$6$};
		\filldraw[fill=white] (N7_2) circle (\RAD) node {$7$};
	\end{tikzpicture}
	
	\caption{Examples of communication protocols in a ring with an odd number of agents ($m=7$): edges with the same colors are activated simultaneously. Three rounds are needed to activate all edges in both protocols.}
	\label{fig:communication_protocol_odd}
\end{figure}

\begin{figure}[b]
	\centering
	\begin{tikzpicture}
		\def\STEP{0.5}
		\def\RAD{0.3}
		
		% Define colors
		\definecolor{blue}{rgb}{0.3010,0.7450,0.9330}
		\definecolor{red}{rgb}{0.9290,0.6940,0.1250}
		
		\coordinate (O) at (0,0);
		
		\coordinate (N1a) at ($(O)+(+\STEP,0)$);
		\coordinate (N1b) at ($(O)+(-\STEP,0)$);
		\coordinate (N2a) at ($(O)+(-5*\STEP,0)$);
		\coordinate (N2b) at ($(O)+(-7*\STEP,0)$);
		\coordinate (N7a) at ($(O)+(+7*\STEP,0)$);
		\coordinate (N7b) at ($(O)+(+5*\STEP,0)$);
		
		% Real links
		\draw[thick,blue] (N1b) -- (N2a);
		\draw[thick,red,style={decorate, decoration={snake, segment length=5, amplitude=2}}] (N1a) -- (N7b);
		
		% Agents
		\draw[dashed] ($0.5*(N1a)+0.5*(N1b)$) ellipse (3.5*\RAD cm and 2*\RAD cm) node[above, yshift=2.5*\RAD cm] {$1$};
		\draw[dashed] ($0.5*(N2a)+0.5*(N2b)$) ellipse (3.5*\RAD cm and 2*\RAD cm) node[above, yshift=2.5*\RAD cm] {$2$};
		\draw[dashed] ($0.5*(N7a)+0.5*(N7b)$) ellipse (3.5*\RAD cm and 2*\RAD cm) node[above, yshift=2.5*\RAD cm] {$7$};
		
		% Masks
		\filldraw[white] ($0.5*(N2a)+0.5*(N2b)+(0,-2.1*\RAD)$) rectangle ($0.5*(N2a)+0.5*(N2b)+(-3.6*\RAD,2.1*\RAD)$);
		\filldraw[white] ($0.5*(N7a)+0.5*(N7b)+(0,-2.1*\RAD)$) rectangle ($0.5*(N7a)+0.5*(N7b)+(3.6*\RAD,2.1*\RAD)$);
		
		% Virtual links
		\draw[thick,black,dotted] (N1a) -- (N1b);
		\draw[thick,black,dotted] (N2a) -- ($(N2b)+(+\RAD,0)$);
		\draw[thick,black,dotted] (N7b) -- ($(N7a)+(-\RAD,0)$);
		
		% Subagents
		\filldraw[fill=white] (N1a) circle (\RAD) node {$1a$};
		\filldraw[fill=white] (N1b) circle (\RAD) node {$1b$};
		\filldraw[fill=white] (N2a) circle (\RAD) node {$2a$};
%		\filldraw[fill=white] (N2b) circle (\RAD) node {$2b$};
%		\filldraw[fill=white] (N7a) circle (\RAD) node {$7a$};
		\filldraw[fill=white] (N7b) circle (\RAD) node {$7b$};
	\end{tikzpicture}
	
	\caption{Example of subagents and virtual communication links. Original agents are represented as dashed ellipsoids and virtual links as black dotted lines.}
	\label{fig:subagents}
\end{figure}
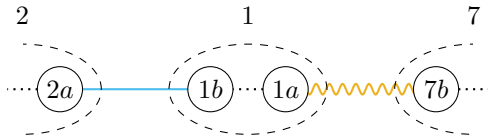

In the case of an odd number $m=2n+1$ of agents, we have to define a proper sequence of mode activations, thus choosing a coloring scheme for the ring. As stated earlier in Section~\ref{sec:setup}, for the odd case we need at least three colors for all edges to be activated in the minimum number of rounds while satisfying the gossip constraint. Whilst in the case of a ring with an even number of agents the coloring scheme is essentially unique (the solutions where two colors are swapped are indeed equivalent, see Figure~\ref{fig:communication_protocol_even}), with three colors we can have multiple coloring configurations. To better clarify this observation, in Figure~\ref{fig:communication_protocol_odd} we report two examples of a coloring scheme with three colors (cyan, yellow and green) for $m = 7$, which translates into two different multigossip sequences for the agents' communication.
Interestingly, the solution that we propose in this paper does not depend on the adopted coloring scheme, which can be arbitrarily chosen but has to be agreed upon execution of the algorithm.

The solution for the case of an odd number of agents is based on the idea of making the number of agents even by considering each agent $i$ as if it were composed by two subagents, say $ia$ and $ib$, which can communicate with the predecessor and successor neighbor of agent $i$ respectively, and are connected together by a \textit{virtual communication link}. The communication link is  virtual because it does not require any actual communication, but just a computation step for agent $i$.
Figure~\ref{fig:subagents} represents a pictorial view of this principle with reference to agent $i=1$ of the ring network at the left panel of Figure~\ref{fig:communication_protocol_odd}. Original agents are represented as dashed ellipsoids and virtual links as black dotted lines.

\begin{algorithm}[t]
	\caption{Algorithm for agent $i$ -- $m$ odd}
	\begin{algorithmic}[1]
		\STATE $x_{ia}(1) \gets$ initial value for agent $i$
		\STATE $x_{ib}(1) \gets$ initial value for agent $i$
		\STATE $k\gets 2$
		\WHILE{$k \leq m$}
			\STATE \textit{\% Exchange information with neighbors}
			\STATE $j \gets \Pre(i)$
			\STATE $x_{ia}(k) \gets (1-\alpha_k) x_{ia}(k-1) + \alpha_k x_{jb}(k-1)$
			\STATE $j \gets \Post(i)$
			\STATE $x_{ib}(k) \gets (1-\alpha_k) x_{ib}(k-1) + \alpha_k x_{ja}(k-1)$
			\STATE $k\gets k+1$
			\STATE \textit{\% Update local estimates}
			\STATE $x_{ia}(k) \gets (1-\alpha_k) x_{ia}(k-1) + \alpha_k x_{ib}(k-1)$
			\STATE $x_{ib}(k) \gets (1-\alpha_k) x_{ib}(k-1) + \alpha_k x_{ia}(k-1)$
			\STATE $k\gets k+1$
		\ENDWHILE
		\RETURN $x_{ia}(m)$
	\end{algorithmic}
	\label{algo:agent_algorithm_odd}
\end{algorithm}

Clearly the topology of this virtual network is still a ring, but it now has an even number $2m$ of agents. Moreover, if we set
\begin{align}\label{eq:init_sub}
x_{ia}(0) = x_{ib}(0) = x_i(0),
\end{align}
then the average of the initial values of the $2m$ subagents equals the average of the initial values of the original $m$ agents. By relabeling the subagents with numbers from $1$ (for subagent 1a) to $2m$ (for subagent $mb$), one can recast this problem to the framework of Section~\ref{sec:solution_even}. Thus, by making each subagent run Algorithm~\ref{algo:agent_algorithm} synchronously with the others, according to Theorem~\ref{thm:correctness_algorithm} we should obtain the average in a finite number of steps $m$, by alternating the edge activation as in Figure~\ref{fig:communication_protocol_even} with $2m$ in place of $m$, with the blue straight links being replaced by the black dotted lines playing the role of the virtual links, and the red wavy links being thought of as the actual communication links, which with reference to Figure \ref{fig:communication_protocol_odd}, would be the cyan, the yellow and the green ones.

Following this principle, however, there is a remaining issue that needs to be addressed, that refers to the fact that simultaneously activating the cyan, yellow and green links would violate the pairwise communication constraint.
This can be alleviated by activating them not all at the same time, but rather following a suitable communication scheme that is compliant with the gossip constraint. For example, we can first activate the cyan straight edges, then the yellow wavy ones, and finally the green spring edges.
This leads to an admissible solution but increases the number of communication rounds, that, unlike the even case, are no longer in correspondence to the algorithm iterations, but rather a multiple of $3$.
In particular, by Theorem~\ref{thm:correctness_algorithm} we know that if the agents in the ring are $2m$, then the average is achieved in $m$ steps. Let $m=2n+1$. If the virtual communication links are activated first, then we have $n+1$ iterations corresponding to virtual communication links and $n$ iterations corresponding to the presence of an actual communication link. However, following the preceding discussion, the number of communication rounds corresponding to the presence of an actual link would be $3n$, i.e., it is a multiple of $3$ of the corresponding number of iterations, to account for the gossip constraint.

\begin{rem}
Due to the fact that $\alpha_1 = 1/2$ (see \eqref{eq:weighting_coeff}), if we start by activating the virtual links, then the first iteration of the algorithm ($k=1$) is not needed since it computes the average of the initial estimates of the two subagents composing each agent which are set equal according to \eqref{eq:init_sub}.
\end{rem}

Algorithm~\ref{algo:agent_algorithm_odd} describes the steps performed by agent $i$. Note that each agent maintains two estimates that are updated based on the information received from its successor and predecessor neighbor and are then combined when the virtual edge is activated.

\section{Conclusions}\label{sec:conclusions}
In this paper we proposed an algorithm for finite time distributed averaging in the case of a ring network of agents, subject to a gossip constraint on communications.
Interestingly, if the number of agents is even, then, consensus to the actual average is achieved in the minimum possible number of iterations, i.e., the diameter of the network, whereas the number of iterations needed in the case where the number of agents os od is higher, but still finite.
Besides being of interest on its own, the proposed algorithm can also be embedded in distributed optimization schemes where computing the average is instrumental to solving the optimization problem (see e.g. \cite{bertsekas1989parallel}) and finite time convergence is then a requirement.
These optimization schemes has been further developed and tailored to large scale systems in new application areas like energy system \cite{dominguez2011distributed}, which could benefit from a distributed finite time averaging algorithm.

\bibliographystyle{IEEEtran}
\bibliography{IEEEabrv,FTring_biblio}

\end{document}